\documentstyle[12pt]{amsart}

\theoremstyle{plain}
\newtheorem{Thm}{Theorem}[section]
\newtheorem{Cor}[Thm]{Corollary}
\newtheorem{Prop}[Thm]{Proposition}
\newtheorem{Lem}[Thm]{Lemma}

\newtheorem{Def}[Thm]{Definition}
\theoremstyle{remark}

\title{On a conjecture of Lange.}
\author{Barbara\ Russo\\ Montserrat\ Teixidor i Bigas}
\address
{Universit\'a degli Studi di Trento\\Dipartimento di Matematica\\
via Sommarive 14\\38050 Povo(TN)\\Italy\\
Mathematics Department\\
Tufts University\\ Medford MA 02155\\ U.S.A.}
\email{russo@@mpim-bonn.mpg.de\\
teixidor@@dpmms.cam.ac.uk}

\begin{document}
\maketitle

\section*{Introduction}
Let $C$ be a projective non-singular curve of genus $g\ge 2$ .
 Let $E$ be a vector bundle of rank $r$ and degree $d$. 
Fix a positive integer $r'<r$. Define 
$$s_{r'}(E)=r'd-r\max _{E'}\{ degE'|rk E'=r', E'\subset E\} $$
Notice that $E$ is stable if and only if $s_{r'}(E)>0$ for every $r'<r$.
 On the other 
hand, for a generic stable $E$
$$r'(r-r')(g-1)\le s_{r'}(E)<r'(r-r')(g-1)+r$$
(cf [L] Satz 2.2 p.452 and [Hi] Th.4.4).
 One can then stratify the moduli space 
$U(r,d)$ of vector bundles of rank $r$ and degree $d$ according
to the value of $s$. Define 
$$U_{r',s}(r,d)=\{ E\in U(r,d)|s_{r'}(E)=s\} $$
We want to study this stratification. A vector bundle $E\in U_{r',s}(r,d)$
can be written in an exact sequence
$$0\rightarrow E'\rightarrow E\rightarrow E''\rightarrow 0$$
with $E',E''$ vector bundles of ranks $r',r''$ and degrees 
$d',d''$ satisfying $r=r'+r'', d=d'+d'', r'd-rd'=r'd''-r''d'=s$.
Note that the condition $s>0$ is equivalent to the inequality of slopes
$\mu (E')<\mu (E'')$. One expects that when this condition 
is satisfied, a generic such extension will yield a stable $E$. 

We call this statement Lange's conjecture (cf. [L]).

The conjecture is now solved and a great deal is known about the geometry 
of the strata : the rank two case is treated in [L,N], the case 
$s\leq \min (r',r'')(g-1)$ in [B,B,R]. In [T1], the result is proved for the 
generic curve and for every curve if $E$ is assumed to be only semistable.
This apparently implies the result also for $E$ stable (cf.[B]).
In [B,L], a proof is provided for $g\ge (r+1)/2$.

The purpose of this paper is to give a simpler proof of the result 
valid without further assumptions. The method is somehow the converse of
the one used by Brambila-Paz and Lange in [B,L]. They start with the 
most general $E$ in $U(r,d)$ and then show that a suitable 
transformation of $E$ gives a new bundle with smaller $s$.
They need to check then that such an $E$ is in fact stable. 
Here, we start with an $E$ with the smallest possible $s$ and 
produce an $E'$ with larger $s$. Stability then comes for free
because $E'$ is more general than $E$. The drawback is that one needs 
to prove existence of stable vector bundles with small $s$ but this
is surprisingly easy.

Our method of proof provides additional information on the 
geometry of the strata. We can prove that $U_{r',s}(r,d)$ is 
contained in the closure of $U_{r',s+r}(r,d)$ as well as the unicity 
of the subbundle (see also [T2])

Our results can be stated in the following

\begin{Thm}
\label{Theorem} 
 Assume that $0<s\le r'(r-r')(g-1) ,s\equiv r'd(r)$.
Write $d'={r'd-s\over r}$. If $g\ge 2$,
 then $U_{r's}(r,d)$ is non-empty, irreducible of dimension 
$$dim U_{r's}(r,d)= r^2(g-1)+1+s-r'(r-r')(g-1)$$
 Moreover, a generic 
$E\in U_{r',s}(r,d)$ can be written in an exact sequence
$$0\rightarrow E'\rightarrow E\rightarrow E''\rightarrow 0$$
with both $E',E''$ stable and $E'$ is the unique subbundle 
of $E$ of  rank $r'$ and degree $d'$.
\end{Thm}

\begin{Thm}
\label{Theorem2} 
If  $s\ge r'(r-r')(g-1)$, every stable vector bundle
has subbundles of rank $r'$ and degree $d'$. Denote by  
$$A_{r',d'}(E)= \{ E'|rk E'=r',deg E'=d',
 E'\subset E ,E' {\rm saturated}\} .$$
Then, for generic  $E$, $A_{r',d'}(E)$ has dimension 
$$dimA_{r',d'}(E) =r'(r-r')(g-1)-s.$$
\end{Thm}
 These results and our methods of proof can be used to study 
twisted Brill-Noether loci. We can show the following

\begin{Thm}
\label{BrillNoether}
 (twisted Brill-Noether for one section).
Let $E$ be a generic vector bundle of rank $r_E$ and degree $d_E$.
Consider the twisted Brill-Noether loci $W^0_{r_F,d_F}(E)$. This is defined 
as the 
 subset of the moduli space $U(r_F,d_F)$ consisting
of those $F$ such that $h^0(F^*\otimes E)\ge 1$.
Then the dimension of $W^0_{r_F,d_F}(E)$ is the expected dimension
given by the Brill-Noether number
$$ \rho ^0_{r_F,d_F}(E)=r_F(r_F-r_E)(g-1)+r_Fd_E-r_Ed_F$$
 if this number is positive and is empty otherwise.
Moreover, when non-empty 
its generic elements considered as maps $F\rightarrow E$
have maximal rank.
\end{Thm}
We also include a proof of Hischowitz's Theorem that states that
the tensor product of two generic vector bundles is non-special.

\bigskip 
Acnowledgments: The first author was partially supported by MURST
  GNSAGA of CNR (Italy), Max-Planck Institut of Bonn and 
 the University of Trento. The second author is visiting the 
Mathematics Department of the University of Cambridge.
 This
collaboration started during the Europroj meeting ``Vector Bundles and 
Equations'' in Madrid. Both authors are members of the Europroj group
VBAC and received support from Europroj
and AGE to attend this conference.

\section{Existence and dimensionality}

In this section we prove the existence of extensions with central
term stable and we compute  the dimension of the set of vector bundles
that fit in such exact sequences. We need several preliminary results.

\begin{Lem}
\label{h^0=0}
 Let $E$ be a stable vector bundle. Assume that we have an 
exact sequence
$$0\rightarrow E'\rightarrow E\rightarrow E''\rightarrow 0.$$
Then $h^0(E^{''*}\otimes E')=0$.
\end{Lem}  
\begin{pf}
A non-zero map $E''\rightarrow E'$ induces an endomorphism of $E$
that is not an homothethy. This is impossible if $E$ is stable. 
\end{pf}

\begin{Thm}[Hirschowitz]
\label{Hirschowitz}
The tensor product of two generic vector bundles is not special.
\end{Thm}
This result was stated and proved in [Hi], 4.6. As this is , 
unfortunately, still unpublished, we provide an alternative proof 
below.
\begin{pf} We shall denote by $r_G,d_G$ the rank and degree of 
a given sheaf say $G$.

By Serre duality, it is enough to show that if $E,F$ are generic
vector bundles, then $h^0(F^*\otimes E)>0$ implies 
$\chi (F^*\otimes E)>0$.

Assume  $h^0(F^*\otimes E)\not= 0$. Then, there is a non-zero 
map $F\rightarrow E$. Denote by $F'$ its kernel, $I$ its image,
 $E''$ its  cokernel. Let $T$ be the torsion subsheaf of $E''$ and 
$\bar E=E''/T$. 
We then have the following exact sequences of sheaves
$$0\rightarrow F'\rightarrow F\rightarrow I\rightarrow 0$$
$$ \begin{array}{ccccccccc}
 & & & & & & 0& & \\
 & & & & & & \downarrow & & \\
 & &0 & & & & T           & & \\
 & &\downarrow & & & & \downarrow & & \\
0&\rightarrow & I&\rightarrow &E&\rightarrow &E''&\rightarrow &0\\
 & &\downarrow & &\downarrow & &\downarrow & & \\
0&\rightarrow &\bar I&\rightarrow &E&\rightarrow &\bar E&\rightarrow &0\\
  & &\downarrow & & & &\downarrow & & \\
 & &T & & & &0 & & \\
 & &\downarrow & &  & & & & \\
 & & 0& & & & & & \\
\end{array}$$
As $T$ is a torsion sheaf, $I$ is determined by $\bar I$, the support 
of $T$ and for every point in the support a map from the fiber 
of $\bar I$ at the point to the basefield.
Hence,
$$dim\{ \bar I \}\le dim \{ I\} +r_I degT$$
As any vector bundle can be deformed to a stable vector bundle, 
(cf.[N,R]Prop.2.6), $F',\bar I, \bar E$ depend at most on 
$r_{F'}^2(g-1)+1, r_{\bar E}^2(g-1)+1, r_{\bar I}^2(g-1)+1$
 moduli respectively.
From \ref{h^0=0} and the stability of $E,F$, 
$h^0(I^*\otimes F)=0, h^0(\bar E^*\otimes \bar I)=0$. 
Notice that $F$ is determined by $F',I$ and an extension class
in $H^1(I^*\otimes F')$ up to homotethy. Similarly, 
$E$ is determined by $\bar I,\bar E$ and an extension class
in $H^1(\bar E^*\otimes \bar I)$  up to homotethy.
From the genericity of the pair $E,F$, we find
$$r_F^2(g-1)+1+r_E^2(g-1)+1=dimU(r_F,d_F)+dimU(r_E,d_E)\le$$
$$\le dimU(r_{F'},d_{F'})+dim U(r_I,d_I)+dim U(r_{\bar E},d_{\bar E})
+r_IdegT+$$
$$+h^1(I^*\otimes F')-1+h^1(\bar E^*\otimes \bar I)-1$$
$$\le (r_{F'}^2+r_I^2+r_{\bar E^2}+r_{F'}r_I+r_Ir_{\bar E})(g-1)
+1+r_{F'}d_I-r_Id_{F'}+r_Id_{\bar E}-r_{\bar E}d_I+r_IdegT$$
This condition can be written as
$$(*)r_{F'}d_I-r_Id_{F'}+r_Id_{\bar E}-r_{\bar E}d_I+r_IdegT-
r_I(r_I+r_{\bar E}+r_{F'})(g-1)-1\ge 0$$
From the genericity of $E,F$ and [L] Satz 2.2, we obtain
$$\mu (I)-\mu(F')\ge g-1, \mu(\bar E)-\mu (\bar I)\ge g-1$$
Hence
$$\mu (\bar E)-\mu (F')\ge 2(g-1)+degT/r_I$$
Equivalently
$$r_{F'}d_{\bar E}-r_{\bar E}d_{F'}\ge 2r_{F'}r_{\bar E}(g-1)+
{r_{F'}r_{\bar E}\over r_I}deg T$$
Adding (*) and this last inequality, we find
$$\chi (F^*\otimes E)\ge 1+r_{F'}r_{\bar E}(g-1)+{r_{F'}r_{\bar E}
\over r_I}deg T+r_{F'}degT\ge 1$$
\end{pf}

\begin{Lem}
\label{cotasubf}
Denote by $V_{r',s}(r,d)$ the set of stable $E$ that can be 
written in an exact sequence of vector bundles
$$0\rightarrow E'\rightarrow E\rightarrow E''\rightarrow 0$$
with $rk E'=r', deg E'=d'$. Assume that $V_{r's}(r,d)$ is non-empty.
 If $s\le r'(r-r')(g-1)$, 
then the generic such $E$ has only a finite number of subbundles
of rank $r'$ and degree $d'$ such that $ r'd-rd'=s$. If $s\ge r'(r-r')(g-1)$, 
then the dimension of the space of subbundles of rank $r'$ and degree
$d'$ of the generic $E$ is at most 
$s-r'(r-r')(g-1)$
\end{Lem}
\begin{pf}
Several proofs of this fact appear in the literature . 
We sketch a proof here for the convenience of the reader.

The set of subbundles of rank $r'$ and degree $d'$ of $E$ is parametrised 
by the quotient scheme of $E$ of the corresponding rank and degree.
The tangent space to this quotient scheme at the point corresponding 
to a bundle $E$ with subbundle $E'$ and quotient $E''$ is 
$H^0(E^{'*}\otimes E'')$. As $E$ is generic, we can assume $E',E''$
generic. Then, from \ref{Hirschowitz}, $E^{'*}\otimes E''$ is non-special.
Hence, if $s\le r'(r-r')(g-1) ,h^0(E^{'*}\otimes E'') =0$ while if
 $s\ge r'(r-r')(g-1)$, then $ h^1(E^{'*}\otimes E'')=0$ and so 
$h^0(E^{'*}\otimes E'')=s-r'(r-r')(g-1)$.
\end{pf}

\begin{Prop}
\label{Irr}
With the notations of \ref{cotasubf}, if $V_{r',s}(r,d)$
 is non-empty, then, it is irreducible
and the generic $E\in V_{r',s}(r,d)$ can be written in an 
exact sequence as above with $E',E''$ stable. Moreover,    
$dimV_{r',s}(r,d)=   \min [r^2(g-1)+1, r^2(g-1)+1+s-r'(r-r')(g-1)]$
 \end{Prop}
\begin{pf}
This proof appears in [T]. We give a sketch here for the convenience
of the reader .

Consider an extension 
$$0\rightarrow E'_0\rightarrow E_0 \rightarrow E''_0\rightarrow 0$$
with $E_0$ stable. From [N,R] Prop.2.6, there are irreducible families of 
vector bundles ${\cal M}', {\cal M}''$ containing $E'_0, E''_0$ respectively
and whose generic member is stable. Consider the universal family
of extensions ${\bf P}$ of an $E''\in {\cal M''}$ by an $E'\in {\cal M'}$.
Consider the open subset $U\subset {\cal M'}\times   {\cal M''}$ consisting 
of those pairs $(E',E'')$ such that $h^0(E^{''*}\otimes E')=0$.
As $\mu (E')<\mu (E'')$, $U$  contains all pairs in which both 
$E',E''$ are stable.
From \ref{h^0=0}, $(E'_0,E''_0)\in U$. As $h^1(E^{''*}\otimes E')$
is constant on $U$, the inverse image ${\bf P}(U)$   of $U$ in
 ${\bf P}$ is irreducible.
 This proves that the given extension can be deformed to an extension 
with both $E',E''$ stable. By the stability of $E_0$, the generic
central term in an extension in ${\bf P}(U)$ is stable. 
Consider the canonical rational map
 $\pi :{\bf P}(U)\rightarrow U(r,d)$. By definition $V_{r',d'}(r,d)$  is the 
image of this map. Hence, it is irreducible. The dimension of ${\bf P}$ 
can be computed as 
$$dim{\bf P}=dim{\cal M}+dim{\cal M}'+h^1(E^{''*}\otimes E')-1=
(r^2-r'r'')(g-1)+1+s$$ 
From \ref{cotasubf}, the fibers of $\pi $ have dimension 
$\max [0,s-r'(r-r')(g-1)]$. Hence, the result follows.  
\end{pf}

\begin{Prop}
\label{Irrsubf}
Assume $s>r'(r-r')(g-1)$ and $E$ is a generic stable vector bundle. 
If $A_{r'd'}(E)$ is non-empty, then it has 
dimension $s-r'(r-r')(g-1)$.
\end{Prop}
\begin{pf}
With the notations in the proof of 
\ref{Irr}, $A_{r',d'}(r,d)$ are the fibers of $\pi$.
Its dimension has been computed already.
\end{pf}  

\begin{Def} 
 Let $E$ be a vector bundle. A vector bundle $\tilde E$ 
is called an elementary transformation  of $E$ if there is an exact 
sequence 
$$0\rightarrow \tilde E \rightarrow E \rightarrow {\bf C}_P \rightarrow 0$$
Here ${\bf C}_P$ denotes the skyscraper sheaf isomorphic to the 
base field with support on the point $P$.

A vector bundle $\bar E$ is called a dual elementary 
transformation of $E$ if $E$ is an elementary
transformation of $\bar E$. Equivalently, the dual of $\bar E$
is an elementary transformation of the dual of $E$ or equivalently
there is an exact sequence
$$\bar E \rightarrow E(Q)\rightarrow {\bf C}_Q^{r-1}\rightarrow 0.$$
\end{Def}

\begin{Lem} 
\label{TE}
 Let 
$$0\rightarrow E' \rightarrow E\rightarrow E''\rightarrow 0$$
be an exact sequence of vector bundles. Then, for a generic elementary
transformation $\tilde E$ of $E$, we have an exact sequence 
of vector bundles 
$$0\rightarrow \tilde E' \rightarrow \tilde E\rightarrow E''\rightarrow 0$$
where $\tilde E'$ is a generic elementary transformation of 
$E$.
\end{Lem} 
\begin{pf} There is an injective map $0\rightarrow E'_P \rightarrow E_P$. 
The elementary transformation depends on the choice of a map
$E_P\rightarrow {\bf C}_P\rightarrow 0$. If this map is generic,
it induces a non-zero map $E'_P\rightarrow {\bf C}_P$.

Hence, we have a diagram
$$\begin{array}{ccccccccc}

 & &0 & &0 & & & & \\
 & &\downarrow & &\downarrow & & & & \\
0&\rightarrow&  \tilde E'&\rightarrow &\tilde E&
\rightarrow &\tilde E''&\rightarrow &0\\
 & &\downarrow & &\downarrow & & \downarrow & & \\
0&\rightarrow&   E'&\rightarrow & E&
\rightarrow & E''&\rightarrow &0\\
 & &\downarrow & &\downarrow & &  & & \\
 & &{\bf C}_P&\rightarrow & {\bf C}_P& & & & \\
 & &\downarrow & &\downarrow & & & & \\
 & & 0         & & 0         & & & & \\
\end{array}$$

This proves the statement.

\end{pf}
\begin{Lem} 
\label{TE*} Let 
$$0\rightarrow E' \rightarrow E\rightarrow E''\rightarrow 0$$
be an exact sequence of vector bundles. Then, for a generic dual
 elementary transformation $\bar E$ of $E$, we have an exact sequence 
of vector bundles 
$$0\rightarrow  E' \rightarrow \bar E\rightarrow \bar E''\rightarrow 0$$
where $\bar E''$ is a generic  dual elementary transformation of 
$E''$.
 \end{Lem}
\begin{pf}:
Dualise the proof above
\end{pf}

\begin{Prop}
\label{exist}
 Let 
$$0\rightarrow E' \rightarrow E\rightarrow E''\rightarrow 0$$
be an exact sequence of vector bundles. Assume that $E$ is stable.
 Then, there exists an exact sequence of vector bundles 
$$0\rightarrow \hat E' \rightarrow \hat E \rightarrow \hat E''\rightarrow 0$$
satisfying 
\begin{description}
\item[i)] $deg \hat E'=deg E'-1, deg \hat E=deg E$\medskip
\item[ii)] $\hat E$ is stable.
\end{description}
\end{Prop} 
\begin{pf}
Take first an elementary transformation
 of the exact sequence
 based at a point $P$. Take next a dual elementary transformation 
based at a point $Q$. From the two Lemmas above, $deg \hat E'=deg E+1$.
We now construct a family of these transformations which contains 
$E$ as one of its members: let the point $Q$ vary until it coincides 
with $P$. Then, with a suitable choice of the dual transformation, 
one can go back to $E$. The existence of this family of 
vector bundles together with the stability of $E$, implies the stability 
of the generic $\hat E$.
\end{pf}
\begin{Cor}
\label{inclusio}
 If $U_{r's}(r,d)$ is non empty,
 then it is contained in the closure of $U_{r', s+r}(r,d)$.
\end{Cor} 
\begin{pf}
Take $E\in U_{r',s}(r,d)$ and consider an exact sequence
$$0\rightarrow E'\rightarrow E\rightarrow E''\rightarrow 0$$
with $E'$ of rank $r'$ and maximal degree $d'$.
In the proof above, we construct a family with special
member $E$ and generic member $\tilde E$ that fits in an exact sequence
$$0\rightarrow \tilde E'\rightarrow \tilde E\rightarrow \tilde E''
\rightarrow 0$$
with $deg(\tilde E')=d'+1$.
Hence, this $\tilde E\in V_{r',s+r}(r,d)$.  
From \ref{Irr}, $V_{r',s+r}(r,d)$    is irreducible and from the
dimensionality statement in \ref{Irr}, $V_{r',s+r}(r,d)\not\subseteq
 V_{r',s-kr}(r,d), k\ge 0$.
   Hence the generic element in $V_{r',s+r}(r,d)$
is in $U_{r',s+r}(r,s)$.
\end{pf}

\begin{Prop}
\label{spetita}
 Let $C$ be a projective non-singular curve of 
genus $g\ge 2$. Consider an exact sequence of vector bundles on $C$
$$0\rightarrow E'\rightarrow E\rightarrow E''\rightarrow 0.$$
Denote by $r',r'',r,d',d'',d$ the ranks and degrees of $E',E'',E$.
 Assume that $E',E''$ are generic stable vector bundles 
of their ranks and degrees.
 If $0<r'd-rd'\le r$, then, the generic such $E$ is stable.
\end{Prop}
\begin{pf} Assume that $E$ is not stable. Let $F$ be a subbundle 
of $E$ such that $\mu (F)\ge \mu (E)$. Up to replacing $F$ by a subbundle
of smaller rank or by its saturation, we can assume $F$ stable and
 $E/F$ without torsion. As $E'$ is stable and
$\mu (E')<\mu (E)$, $F$ gives rise to a non-zero map 
$\phi :F\rightarrow E''$. Denote by $F'$ its kernel, $F''$
 its image.
 Denote by $r_{F'}, r_{F},
r_{F''}, d_{F'}, d_{F}, d_{F''}$ the ranks and degrees of the bundles 
$F',F,F''$.

{\it Claim 1.} $ r_{F''}=r''$.

Proof of  Claim 1: Assume $ r_{F''}<r''$.
 By the genericity of $E''$ this implies 
$r_{F''}d''-r''d_{F''}\ge r_{F''}(r''-r_{F''})(g-1)$ (cf.[L] Satz 2.2).
Equivalently
$$\mu (F'')\le \mu (E'')-(1-(r_{F''}/r''))(g-1)$$
By the initial assumption $r'd-rd'\le r$,
$$(*)\mu (E'')\le \mu (E) +1/r''.$$
As $F$ is a destabilizing subbundle,
$$\mu (E)\le \mu (F)$$
and from the stability of $F$
$$\mu (F) \le \mu (F'')$$
with equality if and only if $F=F''$. 
 Notice that 
$$1/r''-(1-(r_{F''}/r''))(g-1)\le 1/r''-(1-(r_{F''}/r''))\le 0$$
With equalities if and only if $g=2, r_{F''}=r''-1$.
 Puting together   all of the above inequalities, we find that they are all 
equalities. This proves Claim 1 except in the case when all of the following 
properties are satisfied:
\begin{description}
\item[i)] $g=2, r_{F''}=r''-1$
\item[ii)]$\mu(E)=\mu(F)$, $F'=0$ and $F''=F$ is a subsheaf of 
$E''$
\item[iii)] $ (r''-1)d''-r''d_F=r''-1$ 
\end{description}

We shall see at the end of the proof that this situation does not 
correspond to a generic $E$. This will finish the proof of Claim 1.

{\it Claim 2.} $F'=0$

Proof of Claim 2.
Note that $E^*$ satisfies the hypothesis in \ref{spetita} . If $F$ is a
maximal destabilising subbundle of $E$ and we write $G=E/F$,
then $G^*$ is a maximal destabilising subbundle of $E^*$. Then, Claim 2
follows from Claim 1.

From now on, we assume that $F$ is a subbundle of $E''$ and $r_F=r''$

 As $F$ is a destabilising subbundle,  $\mu (F)\ge \mu (E)$.
As $F$ is a subbundle of $E''$,  $\mu (F)\le \mu (E'')$. 
Using (*) and $r_F=r''$, we obtain either $d_F=d''$ or
$d_{F}= d''-1$ and $\mu (F)=\mu (E)=\mu (E'')-1/r''$.
 If $F= E''$,  the sequence splits and the extension
is not generic. If
$d_F=d''-1$, we have an exact sequence
$$0\rightarrow F\rightarrow E''\rightarrow T\rightarrow 0$$
where $T$ is a torsion sheaf of degree one supported at 
one point, say $P$. Then, $F$ is determined by the choice of 
$P$ and a map from $E_P$ to the base field defined up to homothety.
 Therefore the 
number of moduli for such $F$ is at most $r''$.

Consider the pull-back diagram
$$\begin{array}{ccccccccc} 0& \rightarrow &E'&\rightarrow &E&\rightarrow &E''&
\rightarrow &0\cr
 & &\uparrow & &\uparrow & &\uparrow & & \cr
 0& \rightarrow &E'&\rightarrow &E\times _{E''}F&\rightarrow &F&
\rightarrow &0\cr
\end{array} $$

As $F$ is a subsheaf of $E$, the bottom row splits. Hence the top row
corresponds to an element in the kernel of the map
$$ H^1(E^{''*}\otimes E')\rightarrow H^1(F^*\otimes E')\rightarrow 0.$$
This kernel has dimension at most $h^0(T\otimes E')=deg(T)\times rk(E')=r'$.
Therefore the dimension of the subspace of $ H^1(E^{''*}\otimes E')$ that 
may correspond to unstable extensions
is at most 
$$dim\{ F \} +r'\le r''+r'=r.$$ 

On the other hand,
 $$ h^1(E^{''*}\otimes E')=r'd-rd'+r'(r-r')(g-1)=
r'(r-r')(g-1)+r$$
 where the last equality comes from the condition $  \mu (E)=\mu (E'')-1/r''$ .
It is then clear that the generic extension is stable.

We now prove that the special situation at the end of Claim 1 does not occur:
Notice that from condition iii) and the stability of $E''$, $E''$ 
does not have a subbundle of rank $r''-1$ and degree higher than
$d_F$. Hence, we have a pull-back diagram as above but $E/F=L$ is a line
bundle. From condition iii) and \ref{cotasubf},
 there is only a finite number of possible 
$F$ for a given $E''$. From the genericity of $E''$, both $F$ and $L$ are 
generic (and depend only on $E''$ and not on $E'$). Hence, it is enough 
to show that the canonical map
$$ H^1(E^{''*}\otimes E')\rightarrow  H^1(F^*\otimes E')$$
is non-zero. As this map is surjective, this is equivalent to
$ H^1(F^*\otimes E')\not= 0$. From Riemann-Roch
$ h^1(F^*\otimes E')\ge r'(r''-1)+r''(r''-1)(\mu (F)-\mu (E'))>0$
where the last inequality comes from ii) $\mu (E')<\mu (E)$ and
$r_F=r''-1>0$.
\end{pf}

  \begin{pf}. We now prove \ref{Theorem} except for the unicity of the 
subbundle that we
 postpone to next section.
From \ref{cotasubf} , \ref{Irr} and \ref{inclusio} ,
 it is enough to prove the non-emptiness
of $U_{r's}(r,d)$.  Take now any positive $s$. Write $s=ar+\bar s, \  
0<\bar s\le r$. Then, $U_{r',\bar s}(r,d)$ is non empty
 by \ref{spetita}. Applying 
\ref{exist}  $a$-times, we obtain the non-emptiness of $V_{r',s}(r,d)$.
From the definitions of $V_{r',s}(r,d),U_{r',s}(r,d)$, a generic element of
 $V_{r',s}(r,d)$ belongs to an  $U_{r',\tilde s}(r,d)$ for some
$\tilde s\le s$. In order to prove the non-emptiness of 
 $U_{r',s}(r,d)$, it is enough to see that 
  $V_{r',s}(r,d)\not\subset V_{r',\tilde s}(r,d), \tilde s<s$.
 From the dimensionality statement in Lemma 5,
this is true. 
\end{pf}

\begin{pf} The  proof of  \ref{Theorem2} is  similar to the proof
of \ref{Theorem}: use
\ref{spetita}, \ref{exist} and the dimensionality statement in 
\ref{Irr}. 
\end{pf}

\section{Brill-Noether for twisted bundles and unicity of subbundles}
 In this section we prove the result that we stated in the introduction
about twisted  Brill-Noether Theory. The corresponding result 
for the untwisted case (i.e. $E={\cal O}$) is well-known (cf[S] Th IV 2.1).
 We then use this result to show the unicity of the Lange subbundle.
\begin{pf} (of \ref{BrillNoether})
Assume that  $ W^0_{r_F,d_F}(E)$ is non-empty.
Consider an element $F$ in this set. This gives rise to 
a non-zero map $F\rightarrow E$. Denote by $F'$ its kernel, 
$F''$ its image. Then $F''$ is a subbundle of $E$. Moreover,  we have an 
exact sequence
$$0\rightarrow F'\rightarrow F\rightarrow F''\rightarrow 0.$$
Assume first $r_{F''}<r_E$. From \ref{Theorem}, the set of
saturated subbundles 
of $E$ of rank $r_{F''}$ and degree $d_{F''}$ has 
dimension 
$$dim\{ F''\} = r_{F''}(r_{F''}-r_{E})(g-1)+r_{F''}d_{E}-r_{E}d_{F''}$$
if this number is positive and is empty otherwise. The set of 
non-saturated subbundles has dimension smaller than this number.

Consider then the case in which $r_{F''}=r_E$. Then, the quotient 
$E/F''$ is torsion. The choice of $F''$ depends on the choice of the support
of this quotient and for each point $P$ on the support the choice of 
a map (up to homothety) to the base field $E_P\rightarrow {\bf C}$ .
Hence, 
$dim\{ F''\} \le r_E(d_E-d_{F''})$. This coincides with the bound above
for the case $r_{F''}=r_E.$

Any family of vector bundles can be embedded in a family with 
generic member stable(cf [NR] Prop.2.6). Hence, 
$F'$ varies in a parameter space of dimension at most 
$$dim \{ F' \} =r_{F'}^2(g-1)+1.$$
From \ref{h^0=0},  $h^0(F*{''*}\otimes F')=0$. Using Riemann-Roch,
 $$h^1(F^{''*}\otimes F')=
r_{F'}r_{F''}(g-1)+r_{F'}d_{F''}-r_{F''}d_{F'}.$$
The choice of $F$ depends on the choice of the pair $F',F''$ and 
the class of the extension up to scalar. 
Therefore, the dimension of all possible $F$ is bounded by 
$$dim \{ F\} \le [r_{F'}^2+r_{F''}^2-r_{F''}r_E+r_{F'}r_{F''}](g-1)+
r_{F'}d_{F''}-r_{F''}d_{F'}+r_{F''}d_E-r_Ed_{F''}=$$
$$=\rho ^0_{r_F,d_F}(E)-[r_{F'}d_{E''}-d_{F'}r_{E''}-
r_{F'}r_{E''}(g-1)]$$
where $E''$ denotes the quotient of $E$ by $F''$
Let us check that
 $$(*)[r_{F'}d_{E''}-d_{F'}r_{E''}-r_{F'}r_{E''}(g-1)]\ge 0.$$
By the genericity of $E$, if $F''$ exists, then 
$\mu (E'')-\mu (F'')\ge g-1$ (cf Prop. 2.4)
By the stability of $F$, $\mu (F')<\mu (F'')$. Hence, 
$\mu (E'')-\mu (F')\ge (g-1)$. This is equivalent to the inequality
(*) and proves the upper bound for the dimension.
Notice also that (*) vanishes if and only if either $F'=0$
or $F''=E$. In both these cases, the map has maximal rank.

It only remains to prove existence in case the Brill-Noether 
number is positive. If $r_F\not= r_E$, this is equivalent to the existence of 
a stable subbundle or quotient of $E$ and is contained in \ref{Theorem2}.
If $r_F=r_E$, one needs to check that the generic elementary transformation
of a stable vector bundle is again stable. This is well known.
\end{pf}

\begin{Prop}
Let $E$ be a stable vector bundle obtained as a generic extension 
$$0\rightarrow E'\rightarrow E\rightarrow E''\rightarrow 0.$$
Assume that $0<s=r'd-rd'\le r'(r-r')(g-1)$. 
Then the only subbundle of $E$ of rank $r'$ 
and degree $d'$ is $E'$.
\end{Prop}
\begin{pf}
Assume that there were another subbundle $ F'$ of rank $r'$ 
and degree $d'$. Denote by $F''$ the quotient sheaf $E/F'$.
 
 Claim: If $E$ is general, both $F',F''$ are generic vector bundles
 of the given ranks and degrees (i.e. as $E$ varies, $F',F''$ vary
in an open dense subset of the corresponding moduli spaces).

Proof of the claim: If $F''$ had torsion, then $E$ would have a subbundle 
of higher degree and from \ref{Theorem} it could not be general.
If $F'$ or $F''$ were not general or were not stable, then
they would move in  varieties of dimension strictly smaller than 
those parametrising $E',E''$. From \ref{h^0=0}, $h^1(E^{''*}\otimes E')=
h^1(F^{''*}\otimes F')$. From \ref{cotasubf}, every $E$ appears in 
at most a finite number of extensions of an $E''$ by an $E'$.
 Hence, $E$ could not be general. This proves the claim.

We obtain  non-zero maps $ F'\rightarrow E''$ and 
$E'\rightarrow  F''$. From the genericity of 
$E', E''$ and \ref{BrillNoether} , the dimension of the sets of 
these $ F',F''$ is at most
$$dim\{ F'\} =r'(r'-r'')(g-1)+r'd''-r''d'=r'(r'-r'')(g-1)+s$$ 
$$dim \{ F''\} =r''(r''-r')(g-1)+r'd''-r''d'=r''(r''-r')(g-1)+s$$
and both these numbers are positive.
From \ref{h^0=0} and the stability of $E$, $h^0(F^{''*}\otimes F')=0$.
Hence, from Riemmann-Roch
$$h^1(F^{''*}\otimes  F')=r'r''(g-1)-[r''d'-r'd'']=
r'r''(g-1)+s$$
We obtain then a bound for the dimension of the set of $E$ for which 
there is an exact sequence 
$$0\rightarrow  F'\rightarrow E\rightarrow F''\rightarrow 0$$
given by 
$$dim \{ E\} \le (r^{'2}+r^{''2}-r'r'')(g-1)+3s-1$$ 
On the other hand, from \ref{Irr}, 
$$dim \{ E\} =(r^{'2}+r^{''2}+r'r'')(g-1)+s+1$$ 
It follows then that 
$2r'r''(g-1)\le 2s-2$. This contradicts our assumption on $s$.
\end{pf}

\end{document}